\documentclass[a4paper,reqno]{amsart}
\newtheorem*{theorem*}{Theorem}
\newtheorem*{lemma*}{Lemma}
\theoremstyle{remark}
\newtheorem*{remark*}{Remark}

\begin{document}
\author{Peter Zograf}
\address{St.Petersburg Department of the Steklov Mathematical Institute\\
Fontanka 27\\ St. Petersburg 191023 Russia, and
Chebyshev Laboratory of St. Petersburg State University\\
14th Line V.O. 29B\\
St.Petersburg 199178 Russia}
\email{zograf{\char'100}pdmi.ras.ru}

\title{An explicit formula for Witten's 2-correlators}
\thanks{This work was supported by the Russian Science Foundation grant 16-11-10039.}
\begin{abstract}
An explicit closed form expression for 2-correlators of Witten's two dimensional topological gravity is derived in arbitrary genus.
\end{abstract}

\maketitle

Let $\overline{\mathcal{M}}_{g,n}$ be the Deligne-Mumford moduli space of genus $g$ complex stable algebraic curves with $n>0$ distinct marked points. Consider the tautological line bundles $\mathcal{L}_i\to\overline{\mathcal{M}}_{g,n},\; i=1,\ldots,n$. Recall that $\mathcal{L}_i$ is defined fiberwise by $\mathcal{L}_i\left|_{C,x_1,\ldots,x_n}\right.\cong T^*_{x_i}C$, where $C$ is a genus $g$ curve with marked points $x_1,\ldots,x_n$.

Put $\psi_i=c_1(\mathcal{L}_i),\; i=1,\ldots,n,$ and, following Witten \cite{W}, define
\begin{align*}
\langle\tau_{d_1}\ldots\tau_{d_n}\rangle=\int_{\overline{\mathcal{M}}_{g,n}}\psi_1^{d_1}\ldots\psi_n^{d_n}\;,
\end{align*}
where $d_1+\ldots +d_n=3g-3+n$ (we assume that $\langle\tau_{d_1}\ldots\tau_{d_n}\rangle=0$ if any of $d_i<0$). The intersection numbers $\langle\tau_{d_1}\ldots\tau_{d_n}\rangle$ are called correlators of Witten's two dimensional topological gravity. The famous Witten's conjecture \cite{W} (Kontsevich's theorem \cite{K}) claims that the generating function of these numbers (free energy of two dimensional topological gravity) satisfies the KdV (Korteveg-deVries) hierarchy.

Computability of the intersection numbers $\langle\tau_{d_1}\ldots\tau_{d_n}\rangle$ is an important problem. For $g=0$ by a result of Kontsevich \cite{K}
\begin{align*}
\langle\tau_{d_1}\ldots\tau_{d_n}\rangle=\frac{(n-3)!}{d_1!\ldots d_n!}
\end{align*}
(i.~e. a multinomial coefficient). For $g=1$, closed form expressions for the intersection
numbers as sums of the multinomial coefficients were obtained, e.~g., in \cite{KK}. However, no general explicit formulas for the numbers $\langle\tau_{d_1}\ldots\tau_{d_n}\rangle$ are known for $g>1$.

The objective of this note is to derive an explicit closed formula for 2-correlators, i.~e. the numbers $\langle\tau_k\tau_{3g-1-k}\rangle$, for arbitrary 
$g$ (in this case $n=2$ and $k$ ranges from 0 to $3g-1$).
We start with the formula of \cite{LX} just above Proposition 2.6:
\begin{align*}
(2g-1+n)\langle\tau_k\prod_{i=1}^n\tau_{d_i}\rangle&=(2k+3)\langle\tau_0\tau_{k+1}\prod_{i=1}^n \tau_{d_i}\rangle\nonumber\\
&-\frac{1}{6}\langle\tau_0^3\tau_k\prod_{i=1}^n \tau_{d_i}\rangle-\sum_{I\cup J=\{1,\ldots,n\}}\langle\tau_0\tau_k\prod_{i\in I}\tau_{d_i}\rangle\langle\tau_0^2\prod_{j\in J} \tau_{d_j}\rangle\;.
\end{align*}
Adapted for $n=1$ (the case of 2-correlators), it reads
\begin{align}\label{LX}
&2g\langle\tau_k\tau_{3g-1-k}\rangle=\nonumber\\
&(2k+3)\langle\tau_0\tau_{k+1}\tau_{3g-1-k}\rangle
-\frac{1}{6}\langle\tau_0^3\tau_k\tau_{3g-1-k}\rangle-\langle\tau_0\tau_k\rangle\langle\tau_0^2\tau_{3g-1-k}\rangle\;.
\end{align}
We will use the string equation
\begin{align}\label{string}
\langle\tau_0\prod_{i=1}^n \tau_{d_i}\rangle=\sum_{j=1}^n\langle\prod_{i=1}^n \tau_{d_i-\delta_{ij}}\rangle
\end{align}
and the dilaton equation
\begin{align}\label{dilaton}
\langle\tau_1\prod_{i=1}^n \tau_{d_i}\rangle=(2g-2+n)\langle\prod_{i=1}^n \tau_{d_i}\rangle
\end{align}
for Witten's correlators, see \cite{W}.
Applying the string equation \eqref{string} to \eqref{LX}, we easily get that
\begin{align}
(2k+3)&\langle\tau_{k+1}\tau_{3g-2-k}\rangle=(2g-3-2k)\langle\tau_{k}\tau_{3g-1-k}\rangle\nonumber\\
&+\frac{1}{6}\left(\langle\tau_{k-3}\tau_{3g-1-k}\rangle+3\langle\tau_{k-2}\tau_{3g-2-k}\rangle+3\langle\tau_{k-1}\tau_{3g-3-k}\rangle
+\langle\tau_{k}\tau_{3g-4-k}\rangle\right)\nonumber\\
&+\langle\tau_{k-1}\rangle\langle\tau_{3g-3-k}\rangle\;.\label{rec_tau}
\end{align}
From here, using \eqref{string}, \eqref{dilaton} and the fact from \cite{W} that
$$\langle\tau_{3g-2}\rangle=\frac{1}{24^gg!}\;,$$
one can recursively compute 2-correlators $\langle\tau_k\tau_{3g-1-k}\rangle$ for any $g\geq 2$ and $k=0,\ldots,3g-1$. 

For the sake of convenience, let us put
\begin{align}\label{def_a}
a_{g,k}&=\frac{(2k+1)!!(6g-1-2k)!!\langle\tau_{k}\tau_{3g-1-k}\rangle}{(6g-1)!!\langle\tau_{0}\tau_{3g-1}\rangle}\nonumber\\
&=\frac{(2k+1)!!(6g-1-2k)!!}{(6g-1)!!}24^g g!\langle\tau_{k}\tau_{3g-1-k}\rangle\,.
\end{align}
Then in terms of $a_{g,k}$ we can rewrite \eqref{rec_tau} as follows:
\begin{align}
(6g-1-2k)\,a_{g,k+1}=(2g-&3-2k)\,a_{g,k}\nonumber\\
+\frac{4g}{(6g-1)(6g-3)(6g-5)}&\left((2k+1)(2k-1)(2k-3)\,a_{g-1,k-3}\right.\nonumber\\
&+3(2k+1)(2k-1)(6g-1-2k)\,a_{g-1,k-2}\nonumber\\
&+3(2k+1)(6g-1-2k)(6g-3-2k)\,a_{g-1,k-1}\nonumber\\
&\left.+(6g-1-2k)(6g-3-2k)(6g-5-2k)\,a_{g-1,k}\right)\nonumber\\
&\hspace{-1in}+\begin{cases}\frac{g!}{j!(g-j)!}\,\frac{(2k+1)!!(6g-1-2k)!!}{(6g-1)!!}\,,\quad k=3j-1,\vspace{6pt}\\ 0,\quad \text{otherwise}.\end{cases}\label{rec_a}
\end{align}
Using \eqref{string} and \eqref{dilaton}, it is elementary to show that
\begin{align*}
a_{g,0}=1,\quad a_{g,1}=\frac{6g-3}{6g-1}\,.
\end{align*}

Consider now the differences $b_{g,k}=a_{g,k+1}-a_{g,k},\;k=0,\ldots,\left[\frac{3g-1}{2}\right]-1$. Below we derive simple explicit formulas for these numbers.
Actually, we have the following
\begin{lemma*}\label{thm}
The numbers $b_{g,k}$ are explicitly given by the formulas
\begin{align}
b_{g,k}=\frac{(6g-3-2k)!!}{(6g-1)!!}\cdot\begin{cases}
\frac{(6j-1)!!}{j!}\,\frac{(g-1)!}{(g-j)!}\, (g-2j)\,,\quad k=3j-1,\vspace{6pt}\\
-2\,\frac{(6j+1)!!}{j!}\,\frac{(g-1)!}{(g-1-j)!}\,,\quad k=3j,\vspace{6pt}\\
2\,\frac{(6j+3)!!}{j!}\,\frac{(g-1)!}{(g-1-j)!}\,,\quad k=3j+1.\label{diff_a}
\end{cases}
\end{align}
\end{lemma*}
\begin{proof}
Take the difference
\begin{align*}
(6g-3-2k)\,a_{g,k+2}-(6g-1-2k)\,a_{g,k+1}
\end{align*}
and apply formula \eqref{rec_a} to both of its terms. A straightforward computation yields the following recursion for $b_{g,k+1}=a_{g,k+2}-a_{g,k+1}$:
\begin{align}
(6g-3-2k)\,b_{g,k+1}=(2g-&3-2k)\,b_{g,k}\nonumber\\
+\frac{4g}{(6g-1)(6g-3)(6g-5)}&\left((2k+1)(2k-1)(2k-3)\,b_{g-1,k-3}\right.\nonumber\\
&+3(2k+1)(2k-1)(6g-3-2k)\,b_{g-1,k-2}\nonumber\\
&+3(2k+1)(6g-3-2k)(6g-5-2k)\,b_{g-1,k-1}\nonumber\\
&\left.+(6g-3-2k)(6g-5-2k)(6g-7-2k)\,b_{g-1,k}\right)\nonumber\vspace{6pt}\\
&\hspace{-1in}+\begin{cases}\frac{g!}{j!(g-j)!}\,\frac{(2k+3)!!(6g-3-2k)!!}{(6g-1)!!}\,,\quad k=3j-2,\vspace{6pt}\\ 
-\frac{g!}{j!(g-j)!}\,\frac{(2k+1)!!(6g-1-2k)!!}{(6g-1)!!}\,,\quad k=3j-1,\vspace{6pt}\\
0, \qquad k=3j.\end{cases}\label{rec_b}
\end{align}
Now take the values of $b_{g,k}$ given by \eqref{diff_a} and substitute them into \eqref{rec_b}. After a lengthy but elementary computation we see that the numbers 
$b_{g,k}$ satisfy the recursion \eqref{rec_b}. This completes the proof since the recursion \eqref{rec_b} has a unique solution with given initial values.
\end{proof}

By the definition of the numbers $b_{g,k}$ we have 
$$a_{g,k}=a_{g,1}+\sum_{i=1}^{k-1}b_{g,i}\;,\quad k=2,\ldots,\left[\frac{3g-1}{2}\right]-1\;.$$ 
Together with \eqref{def_a} this yields
\begin{theorem*}
The following closed form expression is valid for the 2-correlators: 
\begin{align*}
\langle\tau_{k}\tau_{3g-1-k}\rangle=\frac{(6g-1)!!}{24^g g! (2k+1)!!(6g-1-2k)!!}\left(\frac{6g-3}{6g-1}+\sum_{i=1}^{k-1} b_{g,i}\right)\;.
\end{align*}
Here the numbers $b_{g,i}$ are given by formula \eqref{diff_a}, and $k=2,\ldots,\left[\frac{3g-1}{2}\right]-1$.
\end{theorem*}

\begin{remark*}
A careful analysis of the numbers $b_{g,k}$ performed in \cite{DGZZ}, Sect.~4, implies that for any $g\geq 1$ and $k=2,\ldots,3g-3$ one has
\begin{align*}
\frac{6g-3}{6g-1}<a_{g,k}<1\,.
\end{align*}
\end{remark*}

{\bf Acknowledgements.} The author thanks A.~Zorich for stimulating discussions.

\end{document}